\documentclass[12pt]{article}
\usepackage{amsfonts,amssymb}

\font\scalefont = cmti8
 
\thicklines

\newcount\xmin                
\newcount\xmax
\newcount\ymin
\newcount\ymax

\newcount\gridwidth
\newcount\gridheight

\newcount\nofxpoints          
\newcount\nofypoints          

\newcount\dcnta               
\newcount\dcntb               

\def\begingrid#1#2#3#4{			
	\global\xmin = #1
	\global\ymin = #2
	\global\xmax = #3
	\global\ymax = #4
	\ifnum\xmin > \xmax\errmessage{PATHS: \xmin > \xmax|}\fi
	\ifnum\ymin > \ymax\errmessage{PATHS: \ymin > \ymax|}\fi
	\gridwidth = \xmax
	\gridheight = \ymax
	\advance\gridwidth by -\xmin
	\advance\gridheight by -\ymin
	\nofxpoints = \gridwidth
	\advance\nofxpoints by 1
	\nofypoints = \gridheight
	\advance\nofypoints by 1
	\advance\gridwidth by 4
	\advance\gridheight by 4
	\advance\xmin by -2
	\advance\ymin by -2
	\begin{picture}(\gridwidth, \gridheight)(\xmin, \ymin)
	\advance\gridwidth by -4
	\advance\gridheight by -4
	\advance\xmin by 2
	\advance\ymin by 2
	\global\dcnta = \ymin
	\loop\ifnum\dcnta<\ymax
		\makegridline\dcnta	
		\global\advance\dcnta by 1
	\repeat
	\makegridline\ymax
}
 
\def\makegridline#1{
   \begingroup
   \global\dcntb = \xmin
   \loop\ifnum\dcntb<\xmax
      \put(\dcntb, #1){\circle*{0.1}}
      \global\advance\dcntb by 1
   \repeat
   \put(\xmax, #1){\circle*{0.1}}
   \endgroup
}
 
\def\makexaxis{\thinlines
        \dcnta = \gridwidth
        \advance\dcnta by 1
        \put(\xmin, 0){\vector(1,0){\dcnta}}
        \thicklines
}
\def\makeyaxis{\thinlines
        \dcntb = \gridheight
        \advance\dcntb by 1
        \put(0, \ymin){\vector(0,1){\dcntb}}
        \thicklines
}
\def\makeaxes{\thinlines
        \dcnta = \gridwidth
        \advance\dcnta by 1
        \put(\xmin, 0){\vector(1,0){\dcnta}}
        \dcntb = \gridheight
        \advance\dcntb by 1
        \put(0, \ymin){\vector(0,1){\dcntb}}
        \thicklines
}
 
\def\makexscale{
   \dcnta = \xmin
	\ifnum\dcnta<-9\dcnta=-9\fi
   \loop\ifnum\dcnta<0
		\put(\dcnta,0){\line(0,-1){0.2}}
      \put(\dcnta,-0.5){\makebox(0,0){
		{\scalefont\number\dcnta}}}
      \advance\dcnta by 1
  \repeat
  \dcnta = 1
  \dcntb=\xmax
  \ifnum\dcntb>9\dcntb=9\fi
  \loop\ifnum\dcnta<\dcntb
      \put(\dcnta,0){\line(0,-1){0.2}}
      \put(\dcnta, -0.5){\makebox(0,0){
		{\scalefont\number\dcnta}}}
      \advance\dcnta by 1
  \repeat
  \put(\dcntb,0){\line(0,-1){0.2}}
  \put(\dcntb, -0.4){\makebox(0,0){{\scalefont\number\dcntb}}}
}
\def\makeyscale{
   \dcnta = \ymin
   \loop\ifnum\dcnta<0
      \put(0,\dcnta){\line(-1,0){0.2}}
      \put(-0.4, \dcnta){\makebox(0,0){
		{\scalefont\number\dcnta}}}
      \advance\dcnta by 1
   \repeat
   \dcnta = 1
   \loop\ifnum\dcnta<\ymax
      \put(0,\dcnta){\line(-1,0){0.2}}
      \put(-0.4, \dcnta){\makebox(0,0){
		{\scalefont\number\dcnta}}}
      \advance\dcnta by 1
   \repeat
   \put(-0.4, \ymax){\makebox(0,0){{\scalefont\number\ymax}}}
}

\newsavebox{\skewdowndottedbox}
\newsavebox{\skewupdottedbox}
\savebox{\skewdowndottedbox}(1,1)[lb]{
	\put(0.2,-0.2){\circle*{0.01}}
   \put(0.4,-0.4){\circle*{0.01}}
 	\put(0.6,-0.6){\circle*{0.01}}
	\put(0.8,-0.8){\circle*{0.01}}
}
\savebox{\skewupdottedbox}(1,1)[lb]{
	\put(0.2,0.2){\circle*{0.01}}
   \put(0.4,0.4){\circle*{0.01}}
 	\put(0.6,0.6){\circle*{0.01}}
	\put(0.8,0.8){\circle*{0.01}}
}

\newsavebox{\hplainbox}
\savebox{\hplainbox}{\put(0,0){\line(1,0){1}}}
\newsavebox{\hdottedbox}
\savebox{\hdottedbox}(1,0)[bl]{
	\put(0.2,0){\circle*{0.01}}
	\put(0.4,0){\circle*{0.01}}
	\put(0.6,0){\circle*{0.01}}
	\put(0.8,0){\circle*{0.01}}
}
\newsavebox{\vplainbox}
\savebox{\vplainbox}{\put(0,0){\line(0,1){1}}}
\newsavebox{\vdottedbox}
\savebox{\vdottedbox}(0,1)[bl]{
	\put(0,0.2){\circle*{0.01}}
	\put(0,0.4){\circle*{0.01}}
	\put(0,0.6){\circle*{0.01}}
	\put(0,0.8){\circle*{0.01}}
}
\newsavebox{\vplaindownbox}
\savebox{\vplaindownbox}{\put(0,0){\line(0,-1){1}}}
\newsavebox{\vdotteddownbox}
\savebox{\vdotteddownbox}(0,1)[bl]{
	\put(0,-0.2){\circle*{0.01}}
	\put(0,-0.4){\circle*{0.01}}
	\put(0,-0.6){\circle*{0.01}}
	\put(0,-0.8){\circle*{0.01}}
}
\newsavebox{\skewupbox}
\savebox{\skewupbox}{\put(0,0){\line(1,1){1}}}
\newsavebox{\skewdownbox}
\savebox{\skewdownbox}{\put(0,0){\line(1,-1){1}}}
\newsavebox{\hstep}		
\newsavebox{\vstep}		
\newsavebox{\sstep}		
\newcount\updownincrement

\def\dosteplist{\afterassignment\handlenextstep\let\next=}
\def\handlenextstep{
	\ifx\next\endList
   	\let\next=\relax
	\else
		\ifx\next-
			\put(\dcnta,\dcntb){\usebox{\hstep}}
			\advance\dcnta by 1
		\else
			\ifx\next|
				\put(\dcnta,\dcntb){\usebox{\vstep}}
				\advance\dcntb by\updownincrement
			\else
				\ifx\next/
					\put(\dcnta,\dcntb){\usebox{\sstep}}
					\advance\dcnta by 1
					\advance\dcntb by\updownincrement
				\else
					\errmessage{PATHS: Wrong symbol path argument.}
				\fi
			\fi
		\fi
		\let\next=\dosteplist
	\fi
	\next
}
 
\def\endgrid{\end{picture}}

\newcount\rowcount
\newcount\columncount

\def\begintableau#1#2{									
\global\rowcount=#1
\begin{picture}(#2,#1)(0,0)
\linethickness{0.1pt}
}

\def\endtableau{\end{picture}}

\def\row#1{
	\global\advance\rowcount by -1
	\global\columncount=0
	\dorowlist#1\endList
}

\def\dorowlist{\afterassignment\handlenextentry\let\next=}

\def\handlenextentry{
	\ifx\next\endList
		\let\next=\relax
	\else
		\ifx\next o
			\put(\columncount,\rowcount){\framebox(1,1){\space}}
			\global\advance\columncount by 1
		\else
			\ifx\next-
				\global\advance\columncount by 1
			\else
				\put(\columncount,\rowcount){\framebox(1,1){$\next$}}
				\global\advance\columncount by 1
			\fi
		\fi
		\let\next=\dorowlist
	\fi
	\next
}

\newcount\ca
\newcount\cb

\def\doprintrowlist{\afterassignment\handlenextprintentry\let\next=}

\def\handlenextprintentry{
	\ifx\next\endList
		\let\next=\relax
	\else
		\put(\ca,\cb){\makebox(1,1){$\scriptstyle\next$}}
		\global\advance\ca by 1
		\let\next=\doprintrowlist
	\fi
	\next
}

\newtheorem{num}{\indent\hskip-4.5pt}[section]

\newcommand{\alku}{\bigskip \begin{num}\hskip-6pt .\hskip5pt }
\newcommand{\loppu}{\end{num}}

\newcommand{\be}{\addtocounter{num}{1}\begin{equation}}
\newcommand{\ee}{\end{equation}}

\newcommand{\proof}{{\bf Proof. }\rm }

\font\ff=eusm10 scaled 1200

\begin{document}

\def\esssup{\operatornamewithlimits{ess\,sup}}
\def\Rn{{\mathbb R}^n}
\def\R{{\mathbb R}}
\def\N{{\mathbb N}}
\def\lem{{\rm {\bf Lemma.\hskip 0.5truecm}}}
\def\pro{{\bf Proposition.\hskip 0.5truecm }}
\def\cor{{\bf Corollary.\hskip 0.5truecm}}
\def\thm{{\bf Theorem.\hskip 0.5truecm}}
\def\rem{{\bf Remark.\hskip 0.5truecm}}
\def\define{{\bf Definition.\hskip 0.5truecm}}
\def\lqq{\lq\lq}
\def\rqq{\rq\rq}
\def\hbar{\ \vert \ }
\def\endproof{ $\square$ \bigskip }
\def\period{.\hskip 0.1truecm}
\def\em{\emptyset}

\def\E{\hbox{\ff E}}
\def\k{\hbox{\ff k}}

\pagestyle{myheadings} 
\markright{\sc The electronic journal of combinatorics \textbf{7} (2000), \#Rxx\hfill} 
\thispagestyle{empty}

\title{\bf On Descents in Standard Young Tableaux}
\author{Peter A. H\"ast\"o\footnote{Supported in part by 
the Austrian Academic Exchange Service (\"Osterreichischer akademi- scher
Austauschdienst). 
I'd like to thank C. Krattenthaler for suggesting this topic to me. 
{\hskip 3truecm}
 } 
\\
{\small Department of Mathematics, University of Helsinki,}
\\
{\small P.O. Box 4, 00014, Helsinki, Finland,}
\\
{\small \texttt peter.hasto@helsinki.fi.}}
\date{ \small Submitted: July 9, 2000; Accepted: December 4, 2000}
\maketitle
\abstract{\small In this paper, explicit formulae for the expectation and the
variance of descent functions on random standard Young tableaux are presented. 
Using these, it is shown that the normalized variance, $V/E^2$, is bounded 
if and only if a certain inequality relating the tableau shape to the descent function 
holds.}


\bigskip
\vskip 25pt
\section{Introduction}
\bigskip

In a recent paper, Adin and Roichman defined and studied 
certain descent functions on standard Young tableaux (see \cite{AR}).
They calculated the expectation value
and derived an estimate for the variance of these functions. Their 
results were proved using character theory of symmetric groups.

In this paper, the expectation value and the variance are calculated 
using an elementary method. 
This method is based on the hook-bijection of Novelli, Pak and Stoyanovskii
(cf\period \cite{NPS}). The expressions for the expectation and the variance are 
used to derive a somewhat more precise form of the results in \cite{AR}.

In the following section the necessary definitions are given. Thereafter the 
main results of the paper are stated. In the third section two auxiliary lemmata
are presented. These lead directly to the proofs of the main results and some 
corollaries in the fourth section.

\bigskip
\vskip 25pt
\section{Some definitions and the statement of the main results}
\bigskip

Let $\lambda=(\lambda_1, ..., \lambda_k)$ be a partition of $n$, i.e. a non-increasing sequence of positive integers with sum $n$. We identify
the partition $\lambda$ with its Ferrers diagram (cf\period \cite{ECII}, Section 7.2). A 
{\it standard Young tableau} is a filling of the shape $\lambda$ with the numbers
1, 2, ..., $n$ so that every row and column is increasing. The previous 
statement is to be understood as each number being used exactly once (the \lqq standard\rqq\ part of the name). 
As this paper deals only with standard Young tableaux, it
should be understood that all tableaux (and strings formed from tableaux) 
are made up using each number only once. Figure 1 shows
a standard Young tableau of shape $(4,3,2)$. To every $\lambda$ there corresponds a conjugate partition (cf\period \cite{ECII}, Section 7.2) which will be denoted $\lambda'$. The conjugate partition of
$(4,3,2)$ is $(3,3,2,1)$. 

\begin{figure}[h]
\begin{picture}(15,4)(0,0)
\put(13,0){
	\begintableau{3}{4}
	\row{1346}
	\row{258}
	\row{79}
	\endtableau
}
\end{picture}
\caption{\small A Standard Young Tableau}
\end{figure}

We say that the standard Young tableau $T$ has a {\it descent at} $i$ if 
the entry $i+1$ is 
strictly south (and weakly to the west) of $i$ in $T$. The set of all descents 
in a given tableau $T$ is
denoted by $D(T)$. In the above example $D(T)=\{ 1,4,6,8\}$. 
To every function $f\colon \N \to \R$ there
corresponds a {\it descent function} $d_f(T):= \sum f(i)$ where the sum is
over $D(T)$. Descent functions are generalizations of the classical descent statistics, since $d_f(T)$ equals the descent number of $T$ for $f(n)=1$ and the major 
index of $T$ if $f(n)=n$. 

When talking of {\it random standard Young tableaux} of shape $\lambda$ we 
imply that a uniform distribution is used. The expectation of a descent function over all standard Young tableaux of shape $\lambda$ is denoted by $E_{\lambda}(d_f)$ 
and the variance by $V_{\lambda}(d_f)$.

With these preliminaries we are ready to state our main results:

\alku\label{thm1}\thm For a given $ f\colon \N \to \R$ and $\lambda$ a partition of $n$,
the expectation $E_{\lambda}(d_f) = c_{\lambda'} \sum f(i) $
and the variance 
$$ V_{\lambda}(d_f) = c_{\lambda'} \sum f(i)^2 + 2d_{\lambda'} \sum f(i)f(i+1)$$
$$+ 2(c_{\lambda'}-d_{\lambda}-d_{\lambda'} + e_{\lambda} + e_{\lambda'}) 
\sum_{i-j>1} f(i)f(j) - \left( c_{\lambda'} \sum f(i) \right)^2. $$
Here
$$ N:= n!  \prod_{(i,j)\in \lambda} \{ \lambda_i'+\lambda_j-i-j+1 \}^{-1} $$
is the number of standard Young tableaux of shape $\lambda$ and
$$ c_{\lambda'} := {1 \over 2} \left[ 1 + \sum {\lambda_i'(\lambda_i'-1)
\over n(n-1)} - \sum {\lambda_i(\lambda_i-1)\over n(n-1)}\right]
= \sum_{i\ge j} {\lambda_i'(\lambda_j'-1)\over n(n-1)}, $$
$$ d_{\lambda}:= \sum_{i\ge j\ge k} {\lambda_i(\lambda_j-1)(\lambda_k-2)
\over n(n-1)(n-2)},\ e_{\lambda}:= \sum_{i\ge j\ge k\ge l} 
{\lambda_i(\lambda_j-1)(\lambda_k-2)(\lambda_l-3)\over n(n-1)(n-2)(n-3)},$$
where the sum is over all positive terms, i.e. $\lambda_k>2$ and $\lambda_l>3$.
\loppu 

\alku\rem {\rm
The technique presented in this paper seems to allow the calculation of arbitrary
momenta of $d_f$ (cf. Remark \ref{rem1}). 
However, their expressions will be very complicated, as a comparison
between the expressions for the expectation and the variance leads us to expect.}
\loppu

\alku\label{thm2}\thm Let $f \colon \N \to \R^+$ and let 
$\{\lambda^m \}_{m=1}^{\infty}$ 
be a sequence of partitions.
Then the sequence $ \{V_{\lambda^m}(d_f) / E_{\lambda^m}(d_f)^2 \}_m$ 
is bounded if and only if for all $m$
\be\label{eq1} n\sum_{i=1}^{n-1} f(i)^2 \le c (n-\lambda^m_1) 
\left( \sum_{i=1}^{n-1} f(i)\right)^2, \ee
where $n:= \vert \lambda^m \vert$ and $c$ is a constant independent of $m$.
\loppu

\alku\rem{\rm The reasons for studying $V(X)/E(X)^2$ are two-fold. First, it is the 
variance of the normalized variable $X/E(X)$. Moreover, Chebyshev's inequality is 
used in \cite{AR} to derive a certain type of concentration of $d_f$ if the normalized
variance is bounded. They show that this is the case if $f$ has strictly polynomial 
growth. The results in this paper are more precise, giving the exact condition as to
when the normalized variance is bounded (Theorem \ref{thm2}). However, this 
paper too fails to provide the exact condition under which $d_f$ is concentrated,
since the Chebyshev inequality only provides an upper bound for the distribution. 
}\loppu

\bigskip
\vskip 25pt
\section{Auxiliary results}
\bigskip

The calculation of the expectation is made possible 
by the observation that, for any given $\lambda$, the (average) number of descents
at $i$ is independent of $i$. Similarly, we calculate the variance utilizing 
the fact that the number of co-occurrences of 
descents at $i$ and $i+1$ is independent of $i$ and the number of 
co-occurrences of descents at $i$ and $j$ for $\vert  i - j \vert >1$ is independent of 
$i$ and $j$. The first statement of the following 
lemma is also found in \cite{ECII}, Proposition 7.19.6, 
where it is proved using quasi-symmetric functions.
The lemma follows from Lemma 5.2 of \cite{AR}. However, \cite{AR} uses results from the
character theory of symmetric groups, whereas our proof is wholly elementary. 

\alku\label{lem1}\lem The total number of descents at $i$ ($1\le i<n$) over all 
standard Young tableaux of shape $\lambda$ is independent of $i$. 
The number of co-occurrences of descents at
$i$ ($1\le i<n$) and $j$ ($1\le j<i-1$) is independent of $i$ and $j$. 
The number of co-occurrences of descents at
$i$ and $i+1$ ($1\le i<n-1$) is independent of $i$.
\loppu

\proof Define the partial order $P$ of tableau cells by $c_1 \le c_2$ if
$c_1$ is north-west of $c_2$ (not necessarily strictly). 
Every standard Young tableau corresponds to an extension of the partial order $P$ to
a linear order. Number the cells of
$\lambda$ from left to right, starting from the bottom row (see Figure 2). (This 
labeling is called the superstandard of the shape $\lambda$ in \cite{GR}, p. 219.) Every standard Young tableau also corresponds to a string of numbers formed by reading the numbering in the order determined by the tableau.
(This has been called the {\it inverse reading word} of the tableau.) 
A standard Young tableau has a descent at $i$ exactly when the $i^{\rm th}$
letter of the string 
is greater than the $i+1^{\rm st}$. The standard Young tableau in Figure 2 corresponds to the string $637849152$. This string has descents after the $6$, the $8$, the $9$ and the 
$5$, that is at positions $1$, $4$, $6$ and $8$, respectively. That is to say $D(T)=\{1,4,6,8\}$.

\begin{figure}[h]
\begin{picture}(15,4)(0,0)
\put(8,0){
	\begintableau{3}{4}
	\row{1346}
	\row{258}
	\row{79}
	\endtableau
}
\put(20,0){
	\begintableau{3}{4}
	\row{6789}
	\row{345}
	\row{12}
	\endtableau
}
\end{picture}
\caption{\small A standard Young tableau and the numbering of the cells}
\end{figure}

We arrange the inverse reading words of all tableaux of a given shape on top 
of each other and consider
two adjacent columns of the resulting matrix.
Fix an $i<n-1$. We will prove that the columns $i$ and $i+1$ have the
same number of descents. Let $a<b<c$ be given numbers and consider strings 
corresponding to some standard Young tableau of shape $\lambda$ beginning with some 
string $w$ of length $i-1$ followed by $a$, $b$ and $c$ in some order and ending by some
$v$. We may neglect the strings $wcbav$ (which has descents at both
$i$ and $i+1$) and $wabcv$ (which has neither), should they occur,
since they do not influence the relative number of descents.

Now if $a$, $b$ and $c$ may be chosen freely after $w$, the strings
1) $wcabv$, 2) $wbcav$, 3) $wbacv$ and 4) $wacbv$  will all occur and the 
number of descents due to these is the same in both columns. 
Depending on the shape of the letters of $w$ and $v$ in the tableau 
(for simplicity, we will refer to the cells in which the numbers are as the numbers, so 
the previous means the cells occupied by the elements of the strings $w$ and $v$), 
this may not be
possible. For instance, in the first example in Figure 3, an arbitrary order is possible,
in the second $b$ must precede $a$ allowing the possibilities $wbacv$, $wbcav$ 
and $wcbav$ and so on. Note the that the fifth example is not possible, since
an element of $v$ must precede $b$, which is impossible, since we are considering
strings that end in $v$. Thus not every ordering of $a$, $b$ and $c$ can occur.

\begin{figure}[h]
\begin{picture}(15,4)(0,0)
\put(4,0){
	\begintableau{3}{4}
	\row{wwcv}
	\row{wbv}
	\row{av}
	\endtableau
}
\put(9,0){
	\begintableau{3}{4}
	\row{wwcv}
	\row{bvv}
	\row{av}
	\endtableau
}
\put(14,0){
	\begintableau{3}{4}
	\row{wwcv}
	\row{abv}
	\row{vv}
	\endtableau
}
\put(19,0){
	\begintableau{3}{4}
	\row{wwcv}
	\row{wab}
	\row{wv}
	\endtableau
}
\put(24,0){
	\begintableau{3}{4}
	\row{wcvv}
	\row{wab}
	\row{wv}
	\endtableau
}
\end{picture}
\caption{Examples of constraints (the fifth is impossible)}
\end{figure}
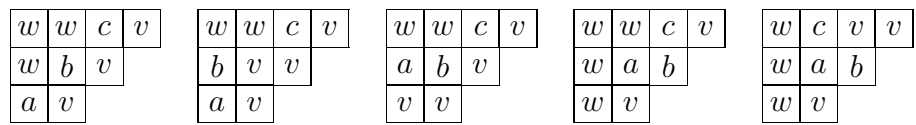

In particular, there are six possibilities with two dependent and one independent cell.
In each exactly two of the four strings above occur and these 
are complementary. For instance, if $c$ must precede $b$, 1) and 4) are
possible. One easily checks the remaining cases as well.
When the order of all three numbers is specified, there are two possibilities.
With a linear order none of the four strings occur. 
In the other case, one number ($b$) precedes or is preceded by the other two. In the 
former case the strings 2) and 3) occur, in the latter 1) and 4), so
again there are equally many descents. These being all the cases, we see
that the columns $i$ and $i+1$ have an equal number of descents, and since 
$i$ was arbitrary the first claim is proved. 

By inspecting the above argument, we see that we have actually proven 
the first co-occurrence claim as well. Fix letters $a<b<c$ and integers $j<i-1<n-2$. 
Forming all possible strings of length $i-1$ and $n-i-2$ without using 
$a$, $b$ and $c$, such that the first has a descent at position $j$ 
and considering these as the $w$ and $v$ of the previous 
paragraphs, we see that the number of descents at $i$ equals the number
of descents at $i+1$ for strings with an $a$, $b$ and $c$ at positions
$i$, $i+1$ and $i+2$ in some order. The first co-occurrence claim follows by 
letting $\{a, b, c\}$ vary over all three-element subsets. 

To prove the last statement, we need to consider strings of the type
$wbcdav$ where $a<b<c<d$ and $w$ has length $i-1$.
There are three interesting 
cases with descents at the first two places: 1a $wdcabv$, 2a $wdbacv$ and 
3a $wcbadv$ as well as three with descents at the last two: 1b $wcdbav$,
2b $wbdcav$ and 3b $wadcbv$. 

There are several cases to analyze. (The analysis that follows is totally elementary,
and the reader may skip to the end of the proof without loss of continuity,
if (s)he is convinced that the proof is \lqq similar\rqq.)
We use the following notation to
describe how $a$, $b$ , $c$ and $d$ relate to the partial order $P$:
$(a,b)$ means that $a$ must precede $b$ in the string and similarly for
more elements. Thus, for instance, $(a,b,d),(c,d)$ means that 
$a$ precedes $b$ which precedes $d$ and $c$ precedes $d$. 
The tables that follow list all possible partial orders on the four elements and which
of the strings 1a-3b from the previous paragraph occur. 

After each table, there are some examples, which show how the first constraint 
of each line might arise. Table 1 lists constraints involving two of the
numbers (cells). There are 12 way of choosing an ordered pair from a 
four-element set, and hence
the table certainly admits all possibilities.

\medskip

\centerline{
\begin{tabular}{cc|cc|cc}
\hbox{} $(a,b)$ & 1a, 3b           & $(a,c)$& 2a, 3b       & $(a,d)$ & 3a, 3b \\
\hbox{} $(b,c)$ & 2a, 2b           & $(b,d)$& 3a, 2b       & $(c,d)$ & 3a, 1b \\
\hbox{} $(b,a)$ & 2a, 3a, 1b, 2b & $(c,a)$& 1a, 3a, 1b, 2b & $(d,a)$ & 1a, 2b, 1b, 2b \\
\hbox{} $(c,b)$ & 1a, 3a, 1b, 3b & $(d,b)$& 1a, 2a, 1b, 3b & $(d,c)$ & 1a, 2a, 2b, 3b \\
\end{tabular}
}
\smallskip
\centerline{Table 1.}

\begin{figure}[h]
\begin{picture}(15,4)(0,0)
\put(7,0){
	\begintableau{4}{4}
	\row{wwwd}
	\row{wwcv}
	\row{abv}
	\row{vv}
	\endtableau
}
\put(12,0){
	\begintableau{4}{4}
	\row{wwwd}
	\row{wbcv}
	\row{avv}
	\row{vv}
	\endtableau
}
\put(17,0){
	\begintableau{4}{4}
	\row{wwwd}
	\row{wwcv}
	\row{bvv}
	\row{av}
	\endtableau
}
\put(22,0){
	\begintableau{4}{4}
	\row{wwwd}
	\row{wwcv}
	\row{wwb}
	\row{av}
	\endtableau
}
\end{picture}
\caption{Examples for Table 1}
\end{figure}
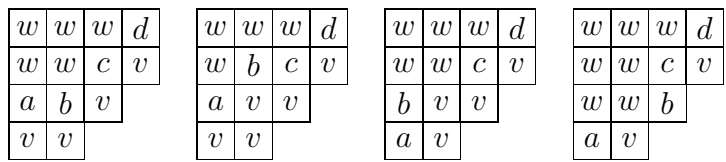

Table 2 contains the relationships between three of the numbers (cells). Since 
constraints with a $v$ cell above a $a$, $b$, $c$ or $d$ cell 
(like the fifth example in Figure 3) and those with a $w$ cell below the same
are impossible, one easily checks that the table lists all cases.

\medskip
\centerline{
\begin{tabular}{cc|cc|cc|cc}
\hbox{} $(a,b,c)$ & $\em$ & $(a,c,d)$& $\em$ & $(a,b,d)$ & $\em$ & $(b,c,d)$& $\em$ \\
\hbox{} $(c,b,a)$ & 3a, 1b & $(d,c,a)$& 1a, 2b & $(d,b,a)$ & 2a, 1a &$(d,c,b)$& 1a, 3b \\
\hbox{} $(a,b),(c,b)$ & 1a, 3b & $(a,b),(d,b)$ & 1a, 3b & 
$(a,c),(d,c)$ & 2a, 3b & $(b,c),(d,c)$ & 2a, 2b \\
\hbox{} $(b,a), (b,c)$ & 2a, 2b & $(b,a), (b,d)$ & 3a, 2b & 
$(c,a), (c,d)$ & 3a, 1b & $(c,b),(c,d)$ & 3a, 1a \\
\end{tabular}
}
\smallskip
\centerline{Table 2.}

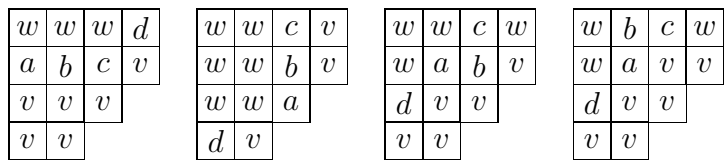
\begin{figure}[h]
\begin{picture}(15,4)(0,0)
\put(7,0){
	\begintableau{4}{4}
	\row{wwwd}
	\row{abcv}
	\row{vvv}
	\row{vv}
	\endtableau
}
\put(12,0){
	\begintableau{4}{4}
	\row{wwcv}
	\row{wwbv}
	\row{wwa}
	\row{dv}
	\endtableau
}
\put(17,0){
	\begintableau{4}{4}
	\row{wwcw}
	\row{wabv}
	\row{dvv}
	\row{vv}
	\endtableau
}
\put(22,0){
	\begintableau{4}{4}
	\row{wbcw}
	\row{wavv}
	\row{dvv}
	\row{vv}
	\endtableau
}
\end{picture}
\caption{Examples for Table 2}
\end{figure}

Tables 3 \& 4 lists the cases where all four of the numbers (cells) are constrained,
specifically Table 4 contains restrictions built up from two pairs and Table 3 
contains the rest. Again, one is easily convinced that these are all the cases. 

\medskip
\centerline{
\begin{tabular}{cc|cc|cc}
\hbox{} $(a,b,c,d)$ & $\em$ & $(d,c,b,a)$ & $\em$ & $(b,c,d),(b,a)$ & $\em$ \\
\hbox{} $(d,c,b),(a,b)$ & 1a, 3b & $(a,b,c),(d,c)$ & $\em$ & $(c,b,a),(c,d)$ & 3a, 1b\\
\hbox{} & & $(c,d,b),(c,a,b)$ & $\em$ \\
\end{tabular}
}
\smallskip
\centerline{Table 3.}

\begin{figure}[h]
\begin{picture}(15,4)(0,0)
\put(9,0){
	\begintableau{3}{4}
	
\row{abcd}
	\row{vvv}
	\row{vv}
	\endtableau
}
\put(14,0){
	\begintableau{3}{4}
	\row{wdvv}
	\row{wcv}
	\row{ab}
	\endtableau
}
\put(19,0){
	\begintableau{3}{4}
	\row{wcdv}
	\row{wab}
	\row{wv}
	\endtableau
}
\end{picture}
\caption{Examples for Table 3}
\end{figure}

\centerline{
\begin{tabular}{cc|cc|cc}
\hbox{} $(a,b),(c,d)$ & $\em$ & $(a,c),(b,d)$& $\em$ & $(a,d),(b,c)$ & $\em$ \\
\hbox{} $(a,b),(d,c)$ & 1a, 3b & $(a,c),(d,b)$& 2a, 3b & $(a,d),(c,b)$ & 3a, 3b \\
\hbox{} $(b,a),(c,d)$ & 3a, 1b & $(c,a),(b,d)$& 3a, 2b & $(d,a),(b,c)$ & 2b, 2b \\
\hbox{} $(b,a),(d,c)$ & 2a, 2b & $(c,a),(d,b)$& 1a, 2b & $(d,a),(c,b)$ & 1a, 1b \\
\end{tabular}
}
\smallskip
\centerline{Table 4.}

\begin{figure}[h]
\begin{picture}(15,4)(0,0)
\put(7,0){
	\begintableau{4}{4}
	\row{wwcd}
	\row{wwvv}
	\row{ab}
	\row{v}
	\endtableau
}
\put(12,0){
	\begintableau{4}{4}
	\row{wwdv}
	\row{wwcv}
	\row{ab}
	\row{v}
	\endtableau
}
\put(17,0){
	\begintableau{4}{4}
	\row{wwcd}
	\row{bvvv}
	\row{av}
	\row{v}
	\endtableau
}
\put(22,0){
	\begintableau{4}{4}
	\row{wwdv}
	\row{wwcv}
	\row{bv}
	\row{a}
	\endtableau
}
\end{picture}
\caption{Examples for Table 4}
\end{figure}

Since for each order-constraint there are
an equal amount of a and b string, we see that, no matter how we choose the cells in 
which to place $a$, $b$, $c$ and $d$,
there are equally many strings with descents at $i$ and $i+1$ as there are strings 
with descents at $i+1$ and $i+2$. The claim now follows by varying $w$ and $v$ 
as above. $\square$

\alku\label{rem1}\rem {\rm
To determine the $n^{\rm th}$ momenta of $d_f$ we need to know the co-occurrence of 
up to $n$ letters. The above method will yield this, however, there will be a lot 
of cases, since we have to consider all the possibilities 
of letters being adjacent or separated by at least one other element.}
\loppu

In order to calculate the expectation and variance of descent functions, we
still need explicit formulae for the invariant numbers put forth in the
previous lemma. To derive these, we count the number of descents at 1, at 
1 and 3, and at 1 and 2. We will use the algorithm from \cite{NPS}
to calculate the relative frequency of standard Young tableaux of these types. 
The algorithm, which was originally devised as a combinatorial proof of the
Stanley hook-length formula, generates standard Young tableaux of uniform distribution.
For ease of reference, the algorithm is described here.

The algorithm of \cite{NPS} starts with a random filling of $\lambda$ with
the numbers 1, ..., $n$. It consists of $n$ steps. In each step there is an
\lqq active\rqq\  element. It is chosen beginning from the rightmost column,
moving up till the column is exhausted, continuing from the bottom of the next
one (to the left) and so on until the upper left corner is reached. 
Having an active element, we compare it with its eastern and 
southern neighbors. If it is the smallest, we move to the next step. Otherwise, 
we exchange the active element with its smaller neighbor and proceed to 
compare it with its new eastern and southern neighbors, continuing the exchanging process 
till the active element is the smallest of the
three numbers. Then we move to the next step.

The algorithm obviously ends in a standard Young tableau, and 
\cite{NPS} tells us that every standard Young tableaux of shape $\lambda$ will be generated exactly $n!/N$ times ($N$  stands for the total number of standard Young tableaux of 
shape $\lambda$ and is given explicitly in the following lemma).
By means of this algorithm we derive:

\alku\label{lem2}\lem For $1\le i<n$, the total number of descents at $i$ 
over all standard Young tableaux of shape $\lambda$ is $Nc_{\lambda'}$. 
The number of co-occurrences of descents at $i$ ($1\le i<n$) and $j$ ($1\le j<i-1$) 
equals $ N(c_{\lambda'}-d_{\lambda}-d_{\lambda'} + 
e_{\lambda} + e_{\lambda'})$.
The number of co-occurrences of descents at
$i$ and $i+1$ ($1\le i<n-1$) is $ Nd_{\lambda'}$ where
$$ c_{\lambda'} := {1\over 2} \left[ 1 + \sum {\lambda_i'(\lambda_i'-1)
\over n(n-1)} - \sum {\lambda_i(\lambda_i-1)\over n(n-1)}\right]
= \sum_{i\ge j} {\lambda_i'(\lambda_j'-1)\over n(n-1)}, $$
$$ d_{\lambda}:= \sum_{i\ge j\ge k} {\lambda_i(\lambda_j-1)(\lambda_k-2)
\over n(n-1)(n-2)},\ e_{\lambda}:= \sum_{i\ge j\ge k\ge l} 
{\lambda_i(\lambda_j-1)(\lambda_k-2)(\lambda_l-3)\over n(n-1)(n-2)(n-3)}$$
and 
$$ N:= n! \prod_{(i,j)\in \lambda} \{ \lambda_i'+\lambda_j-i-j+1 \}^{-1} $$
is the number of standard Young tableaux of shape $\lambda$.

\loppu

\proof Let us inspect the situation for descents of 1. Choose two cells
of $\lambda$. These cells will hold the numbers 1 and 2. 
We generate the rest of the numbers randomly (of uniform distribution). Then 
we use the \cite{NPS} algorithm on this numbering. 

We will discern three cases:

1) Both chosen cells are in the same column. In this case, there will be a descent 
no matter how we place 1 and 2 in the cells. 

2) One cell is in the top row and the other strictly east from it. In 
this case there will never be a descent.

3) Otherwise there will be a comparison of some active element $x$
with both $1$ and $2$ (Figure 8) during the algorithm. Whether there is
a descent on 1 depends on how 1 and 2 are placed in the chosen cells. 

\begin{figure}[h]
\begin{picture}(15,4)(0,0)
\put(8,0){
	\begintableau{2}{2}
	\row{x1}
	\row{2}
	\endtableau
}
\put(11,1){$\to$}
\put(12,0){
	\begintableau{2}{2}
	\row{1x}
	\row{2}
	\endtableau
}
\put(15,1){,}
\put(16,0){
	\begintableau{2}{2}
	\row{x2}
	\row{1}
	\endtableau
}
\put(19,1){$\to$}
\put(20,0){
	\begintableau{2}{2}
	\row{12}
	\row{x}
	\endtableau
}
\end{picture}
\caption{\small The third case}
\end{figure}

The relative frequency of 1) and 2) is respectively 
$\sum \lambda_i'(\lambda_i'-1) /(n(n-1))$ and 
$\sum \lambda_i(\lambda_i-1) /(n(n-1))$. To calculate the relative frequency, 
$c_{\lambda'}$, we take $1/2$ and correct it by adding half of 1) and subtracting
half of 2). The second expression for $c_{\lambda}'$ follows from an argument similar to
that which follows.

For the third claim (the second follows, below) 
we count the number of tableaux with descents
at $1$ and $2$. It is easy to see that the initial situations
that will lead to these descents are those with 1 (weakly) to the right of 2 
which is to the right of 3, 2 is not in the uppermost row
and 3 not in the two uppermost rows. The number of possible such combinations is
$$ \sum_{i\ge j\ge k} \lambda_i'(\lambda_j'-1)(\lambda_k'-2)(n-3)!,$$
which divided by $n!$ yields the fraction
of tableaux of shape $\lambda$ with descents at both 1 and 2. 

We see that $d_{\lambda}$, $e_{\lambda}$ and $e_{\lambda'}$ enumerate
the tableaux with 1, 2 and 3 in the first row, 1, 2, 3 and 
4 in the first row and 1, 2, 3 and 4 in the first column, respectively. Now, the number
of co-occurrences of standard Young tableaux with descents at 1 and at 3 is the number of standard Young tableaux with descents at 1
less the number of standard Young tableaux with descent at 1 but without descent at 3. 
The latter number equals tableaux beginning with 1, 4 / 2 / 3 (meaning 1 and
4 in the first, 2 in the second and 3 in the third row) or with
1, 3, 4 / 2. But these equal respectively the number of tableaux 
with 1 / 2 / 3 less those with 1 / 2 / 3 / 4 (see Figure 9) and
those beginning 1, 2, 3 less those beginning with 1, 2, 3, 4, from which 
the claim follows. \endproof

\begin{figure}[h]
\begin{picture}(15,4)(0,0)
\put(8.5,2){$\#$}
\put(9,0){
	\begintableau{4}{3}
	\row{14o}
	\row{2o}
	\row{3o}
	\row{o}
	\endtableau
}
\put(12.7,2){$= \#$}
\put(14,0){
	\begintableau{4}{3}
	\row{1oo}
	\row{2o}
	\row{3o}
	\row{o}
	\endtableau
}
\put(17.7,2){$-\ \#$}
\put(19,0){
	\begintableau{4}{3}
	\row{1oo}
	\row{2o}
	\row{3o}
	\row{4}
	\endtableau
}
\end{picture}
\caption{See text for details}
\end{figure}

\bigskip
\vskip 25pt
\section{The proofs of the main results}
\bigskip

{ \bf Proof of Theorem 2.1} 
{\rm The claim concerning the expectation follows directly from Lemmata \ref{lem1} and \ref{lem2}. The variance is expressed in
the standard form $V(X)=E(X^2)-E(X)^2$. Since only the co-occurrence of the numbers
matters, this also follows easily from the previous lemmata.} \endproof

{ \bf Proof of Theorem 2.3}
{\rm For convenience, we will denote $\vert \lambda^m \vert$ by $n$ and omit the 
superscript from $\lambda^m$ and write $\lambda_1$ for its first component when 
there is no danger of confusion. We will assume without loss of generality 
that $n\ge 10$.

We start by deriving upper and lower bounds for $c_{\lambda'}$.
Let $q:=\lambda_1/n$ and assume $q<1$ (the claim is trivial otherwise). If $q\ge1/2$ then 
\be\label{c_lowbound1} c_{\lambda'}\ge {1 \over 2} - {q(nq-1)\over 2(n-1)} -
{(1-q)((1-q)n-1) \over 2(n-1)} = { n \over (n-1)}q(1-q) > q(1-q). \ee
If $\lambda_i\le l\le n/2$ and $\Sigma \lambda_i =n$ then 
$$\sum \lambda_i^2\le l+ ... +l + (n-pl) + 0 + ... + 0 \le n^2/2,$$
where the number of $l$-terms is $p$ and $0\le n-pl < l$. This is easily
established by taking some other configuration and increasing the largest term $<l$ by 1
and decreasing the smallest term $>0$ by one. Now the sum is still $n$, but the 
sum of the squares is larger. Since this is a finite process, the maximum 
is as indicated. It follows that if $q<1/2$, 
then $c_{\lambda'}$ will be larger than the least $c_{\lambda'}$ of 
tableaux with $q\le 1/2$, i.e. $c_{\lambda'}\ge 1/4$. Therefore, 
\be\label{c_lowbound} c_{\lambda'}\ge \min \{1/4, (1-q)/2\}. \ee
The upper bound is simpler: 
\be\label{c_highbound} c_{\lambda'} \le {1 \over 2} + { (1-q)((1-q)n+1)\over 2(n-1)} - 
{q(nq-1)\over 2(n-1)} = \ee
$$ ={(1+ (1-q)^2 - q^2)n + 1-q-1+q \over 2(n-1)} = {n \over n-1}(1-q) \le 10(1-q)/9.$$

According to Theorem {\ref{thm1}} we have three terms to bound (the last 
one is already constant, $-1$). We will first show that the third term (corresponding
to descents at $i$ and $j$ for $\vert i - j \vert >1$) is bounded, which is 
equivalent to showing that
$$ (c_{\lambda'}-d_{\lambda}-d_{\lambda'} + e_{\lambda} + e_{\lambda'})
\sum_{i-j>1} f(i)f(j) \le k c_{\lambda'}^2 \left( \sum_{i=1}^{n-1} f(i)\right)^2. $$
Since $f>0$, the sum on the left is less than the square of the 
sum on the right, and we may 
disregard these sums. Since $d_{\lambda'}> e_{\lambda'}$ it suffices to show that there
exists a $k$ independent of $\lambda$ such that
\be\label{former} c_{\lambda'}-d_{\lambda} + e_{\lambda}  \le k c_{\lambda'}^2\ee
holds.

If $q\le 1/2$, then $c_{\lambda'} \ge 1/4$ and (\ref{former}) holds with $k=4$, 
since $d_{\lambda}> e_{\lambda}$. 
Thus we may assume without loss of generality that $q > 1/2$. Let 
$\xi := c_{\lambda'}/(1-q)$. It follows from (\ref{c_lowbound1})
and (\ref{c_highbound}) that $q < \xi < n/(n-1)$. Write $d_{\lambda}$ as
$$ d_{\lambda} = q_1q_2+ q_2\sum_{i\ge j\ge 2} {\lambda_i(\lambda_j-1)\over n(n-1)} + 
\sum_{i\ge j\ge k\ge 2} {\lambda_i(\lambda_j-1)(\lambda_k-2)\over n(n-1)(n-2)}$$
and 
$$ e_{\lambda}= q_1q_2q_3+ q_2q_3\sum_{i\ge j\ge 2} {\lambda_i(\lambda_j-1)
\over n(n-1)} + q_3\sum_{i\ge j\ge k\ge 2} {\lambda_i(\lambda_j-1)(\lambda_k-2)\over n(n-1)(n-2)}+ $$
$$\sum_{i\ge j\ge k\ge l\ge 2} 
{\lambda_i(\lambda_j-1)(\lambda_k-2)(\lambda_l-3)\over n(n-1)(n-2)(n-3)},$$
where $q_i := (nq-i)/(n-i)$ for $i=1,2,3$.
Let $\mu$ be the shape obtained from $\lambda$ by removing the first row. 
Then the sums in the previous equations correspond to $c_{\mu}$, $d_{\mu}$
and $e_{\mu}$. Set
$r:=\vert \mu \vert = n - \lambda_1=n(1-q)$. Combining these equations, we see that
(\ref{former}) is equivalent to
$$ \xi (1-q) + (q_3-1) \left( q_1q_2 + q_2 {r(r-1)\over n(n-1)} c_{\mu} + 
{r(r-1)(r-2)\over n(n-1)(n-2) } d_{\mu} \right) + $$
$$ {r(r-1)(r-2)(r-3)\over n(n-1)(n-2)(n-3) } e_{\mu} \le k {\xi}^2 (1-q)^2.$$
Dividing through by $1-q$, we see, noting that $1-q=r/n$ and that $(q_3-1)/(1-q)=-n
/(n-3)$, that (\ref{former}) is equivalent to
$$ (\xi - n/(n-3)) + {n\over n-3} \left( 1- q_1q_2 - q_2 {r(r-1)\over n(n-1)} c_{\mu} - 
{r(r-1)(r-2)\over n(n-1)(n-2) } d_{\mu} \right) + $$
$$ {(r-1)(r-2)(r-3)\over (n-1)(n-2)(n-3) } e_{\mu} \le k {\xi}^2 (1-q).$$
Since $\xi - n/(n-3) < n/(n-1)-n/(n-3)<0$, we may drop $\xi - n/(n-3)$. After that, we divide by $1-q$ again:
$$ {n\over n-3} \left( {(1+q)n^2-3n \over (n-1)(n-2)}- q_2 {r-1\over n-1} c_{\mu} - 
{(r-1)(r-2)\over (n-1)(n-2)} d_{\mu} \right) + $$
$$ {n(r-1)(r-2)(r-3)\over r(n-1)(n-2)(n-3) } e_{\mu} \le k {\xi}^2 .$$
Since $e_{\mu}$ denotes a fraction of tableaux with 
a specific pattern of descents, it is less than one.
Since $ n(r-1)\le r(n-1)$ and $1+q\le 2$, 
the previous inequality is implied by
$$ {(2n-3)n^2 \over (n-1)(n-2)(n-3)} + {(r-2)(r-3)\over (n-2)(n-3) } \le k q^2,$$
which certainly holds for $k=20$, since we assumed $n\ge 10$ and $q>1/2$. 

The second term in $ V_{\lambda}(d_f) / E_{\lambda}(d_f)^2 $ is dominated by twice the 
first. Hence all that remains to do is to bound the first term, that is, to show that:
\be\label{eq2} \sum_{i=1}^{n-1} f(i)^2 \le k c_{\lambda'} \left( \sum_{i=1}^{n-1} f(i)\right)^2, \ee
where $k$ is some constant independent of $\lambda$ and $f$. 
By ({\ref{eq1}}) and $c_{\lambda'} \ge (1-q)/4$ we get 
$$ \sum_{i=1}^{n-1} f(i)^2 \le c (1-q) \left( \sum_{i=1}^{n-1} f(i)\right)^2\le
4c c_{\lambda'} \left(\sum_{i=1}^{n-1} f(i)\right)^2,$$
so (\ref{eq2}) holds with $k=4c$, where $c$ comes from (\ref{eq1}).

For the necessity of condition (\ref{eq1}) we assume that the normalized variance is 
bounded. Again by Theorem \ref{thm1} this entail that three positive terms are 
bounded, hence each of the terms is bounded. 
To say that the first one is bounded is equivalent to asserting that (\ref{eq2}) holds.
Then (\ref{eq1}) follows if we show that $c_{\lambda'}  \le k'(1-q)$. But this
follows from (\ref{c_highbound}), so we are done. $\square$}
 
\alku\cor\label{cor3} Let $\{ \lambda^m \}_m$ and $f$ be as in Theorem {\ref{thm2}}.
The normalized variance $ V_{\lambda}(d_f) / E_{\lambda}(d_f)^2 $ 
is bounded for all sequences $\{ \lambda^m \}_m$ if and only if 
$$ n\sum_{i=1}^{n-1} f(i)^2 \le c \left( \sum_{i=1}^{n-1} f(i)\right)^2.$$
or equivalently,
$$ \Vert f \Vert_2 \le c \Vert f \Vert_1.\ \square$$
\loppu

\end{document}